%% file: main.tex
\documentclass[reqno, oneside]{amsart}
\usepackage{etex}
\usepackage{chemarr}
\usepackage{hyperref}
\usepackage{comment}

\usepackage[a4paper]{geometry}

\usepackage{latexsym, amscd, amsfonts, eucal, mathrsfs, amsmath, 
amssymb, mathabx, amsthm, xr, stmaryrd, color, enumerate, tikz}
\usepackage{appendix}
\usepackage{multicol}
\usepackage[all]{xy}
\usepackage{hyperref}
\usepackage{fullpage}

\theoremstyle{definition}

\theoremstyle{definition}

\newtheorem*{rmk}{Remark}

\newtheorem*{dfn*}{Definition}
\newtheorem*{axm*}{Axiom}
\newtheorem*{ntn*}{Notation}
\newtheorem*{exm*}{Example}
\newtheorem*{exr*}{Exercise}
\newtheorem*{int*}{Intuition}
\newtheorem*{qst*}{Question}
\newtheorem*{rmk*}{Remark}

\theoremstyle{plain}

\newtheorem{thm}{Theorem}

\newtheorem{lem}[thm]{Lemma}

\newtheorem*{thm*}{Theorem}
\newtheorem*{prop*}{Proposition}
\newtheorem*{cor*}{Corollary}
\newtheorem*{lem*}{Lemma}
\newtheorem*{cnj*}{Conjecture}

\usepackage{tikz}
\usetikzlibrary{matrix,arrows,decorations}
\usepackage{tikz-cd}

\usepackage{adjustbox}

\begin{document}

\title{Quotients of even rings}
\author{Jeremy Hahn and Dylan Wilson}

\begin{abstract} We prove that if $R$ is an $\mathbb{E}_2$-ring
with homotopy concentrated in even degrees, and $\{x_j\}$ is any
sequence of elements in $\pi_{2*}(R)$, then $R/(x_1,x_2,\cdots)$ admits the structure
of an $\mathbb{E}_1$-$R$-algebra.  This removes an assumption, common in the literature, that $\{x_j\}$ be a regular sequence.
\end{abstract}

%----------------------------------------------------------------------%
%----------------------------------------------------------------------%

\maketitle

\vbadness 5000

%----------------------------------------------------------------------%

\input{introduction.tex}

\input{reduction.tex}

\input{augmentation.tex}

\bibliographystyle{amsalpha}
\nocite{*}
\bibliography{Bibliography}

\end{document}

%% file: introduction.tex
\section{Introduction}
Let $R$ be a homotopy commutative ring and $x \in\pi_*R$ an element. 
Then one can ask whether the $R$-module cofiber $R/x$ admits a product compatible with the
action of $R$, up to homotopy. Of particular interest has been
the case when $R$ has homotopy groups concentrated in even degrees.
Variations and special cases of this question
have been considered throughout the years:

\begin{itemize}
\item Bass and Sullivan \cite{bassul} realized various quotients of the complex
bordism ring
$\mathrm{MU}_*$ as coming from homology theories
via a geometric construction. Placing products on the resulting homology
theories was a delicate process, studied early in various references,
such as \cite{shimyag, mor, mir}.
\item Elmendorf-Kriz-Mandell-May \cite{ekmm}, as an application of their new,
highly structured models for ring spectra, showed that various quotients
of $\mathrm{MU}$ could be realized as having multiplications,
often associative or commutative, inside the homotopy category
$h\mathsf{Mod}_{\mathrm{MU}}$.
\item Strickland \cite{strick} generalized the results of \cite{ekmm} by proving,
among other things, that if $R$ is an even, $\mathbb{E}_{\infty}$-ring spectrum,
and $x \in \pi_*R$ is not a zero divisor, then $R/x$ admits the structure of
an associative ring in $h\mathsf{Mod}_R$.
\item Angeltveit \cite{angel} showed that if $R$ is an $\mathbb{E}_{\infty}$-ring
and $x \in \pi_* R$ is not a zero-divisor, then $R/x$ further admits the structure of an $\mathbb{E}_1$-$R$-algebra.  He concludes that, if $(x_1, x_2, ...)$ is a regular sequence in $\pi_*R$, then
$R/(x_1, x_2, ...)$ admits the structure of an $\mathbb{E}_1$-$R$-algebra.
We note that his proof works just as well in the case when
$R$ is an $\mathbb{E}_2$-ring.

%The purpose of this brief note is to remove the assumption, highlighted in the works of Strickland and Angeltveit, that $x$ be a non-zero divisor, or that $(x_1,x_2,\cdots) \in \pi_*(R)$ be a regular sequence.

\item If $x$ is a zero-divisor, then $\pi_*(R/x)$ is not $\pi_*(R)/x$, because the relevant long-exact sequence is not short-exact.  In particular, $\pi_*(R/x)$ will not be concentrated in even degrees, and this significantly interferes with any direct attempt to generalize the arguments of Strickland and Angeltveit.  Nonetheless, the theory of Thom spectra has been used to remove the assumption in certain special cases.  The strongest result of this form is recent work of Basu-Sagave-Schlichtkrull \cite{bss}, who give a direct construction of $R$-algebra quotients $R/(x_1, x_2, ...,x_n)$
in the case when $R$ is an even, $\mathbb{E}_{\infty}$-ring and $|x_i| = 2i$,
with no assumption that the sequence be regular.
We remark that their proof works also
in the case when $R$ is an $\mathbb{E}_2$-algebra, 
using, for example, the machinery
in \cite{omar-tobias}.
\end{itemize}

We offer the following addition to the above list, which is a simultaneous strengthening of Angeltveit's result and the result of Basu-Sagave-Schlichtkrull. 

\begin{thm}\label{thm:main} Let $R$ be an $\mathbb{E}_2$-ring with 
$\pi_{2k-1}R=0$ for all $k \in \mathbb{Z}$.
Let $J=\{x_j\}$ be any sequence
of elements in $\pi_{2*}R$. Define
	\[
	R/J:= \bigotimes_{R} R/x_j.
	\]
Then $R/J$ admits the structure of an $\mathbb{E}_1$-algebra in
left $R$-modules. 
\end{thm}

We will reduce our theorem to Angeltveit's by way of the following:

\begin{lem}\label{thm:augment} For fixed $k \in \mathbb{Z}$, let 
$$A=S^0[S^{2k}]$$ denote the free $\mathbb{E}_1$-ring spectrum generated by $S^{2k}$. Then $A$ admits
the structure of a nonnegatively graded $\mathbb{E}_2$-ring, with degree $0$ component $A_0 = S^0$.  In particular, there is an $\mathbb{E}_2$-algebra augmentation
$$A \longrightarrow S^0.$$
\end{lem}

\begin{rmk} We imagine that Lemma \ref{thm:augment} is well-known, but could not find the statement in the literature so we provide a proof below.
\end{rmk}

\begin{rmk}
In the case that the sequence of elements $\{x_j\}$ appear in nonnegative degree, it is possible to use Thom spectra, elaborating on the argument given in \cite{bss}, to give a proof of Theorem \ref{thm:main}.  However, because $\text{gl}_1(R)$ depends only on the connective cover of $R$, it seems impossible for such techniques to handle the case that the $x_j$ appear in negative degrees, and this case is essential to forthcoming applications of the authors.
\end{rmk}

%% file: reduction.tex
\section{Reduction}

We defer the proof of Lemma \ref{thm:augment} to the next section
and explain how to deduce Theorem \ref{thm:main}.

With notation as in the statement of the main theorem, define
	\[
	S^0[J] := \bigwedge_{j} S^0[S^{|x_j|}].
	\]
Denote by $t_j$ the class $S^{|x_j|} \to S^0[J]$.

By Lemma \ref{thm:augment}, this is an augmented $\mathbb{E}_2$-ring.
Thus
	\[
	R[J] := R \wedge S^0[J]
	\]
is an $\mathbb{E}_2$-ring equipped with an $\mathbb{E}_2$-map
to $R$. (Beware, however, that we cannot even make sense of
$R[J]$ being an $\mathbb{E}_2$-$R$-algebra in general, since
$\mathsf{LMod}_R$ is only guaranteed to be $\mathbb{E}_1$-monoidal).

Now note that $R[J]$ is an $\mathbb{E}_2$-ring with homotopy concentrated
in even degrees and that
	\[
	(x_1-t_1, x_2-t_2, ...) \subseteq \pi_*R[J]
	\]
is an ideal generated by a regular sequence. Thus, by Angeltveit's theorem, we may promote
	\[
	R[J]/(x_1-t_1, x_2-t_2, ...)
	\]
to an $\mathbb{E}_1$-$R[J]$-algebra. Since the augmentation $R[J] \to R$
is an $\mathbb{E}_2$-map,
base change is $\mathbb{E}_1$-monoidal and we deduce
that
	\[
	R \otimes_{R[J]} \left(R[J]/(x_1-t_1, x_2-t_2, ...)\right)
	\]
becomes an $\mathbb{E}_1$-$R$-algebra. Rearranging the tensor product
yields Theorem \ref{thm:main}.

%% file: augmentation.tex
\section{Augmentation}

We are left with proving Lemma \ref{thm:augment}, which states that
$S^0[S^{2k}]$ admits the structure of a nonnegatively graded
$\mathbb{E}_2$-ring with $S^0$ in degree zero. 

Recall that the homotopy
theory
of graded spectra is given by 
$\mathsf{Fun}(\mathbb{Z}^{\mathrm{ds}}_{\ge 0}, \mathsf{Sp})$ where
$\mathbb{Z}_{\ge 0}^{\mathrm{ds}}$ denotes the set of nonnegative
integers regarded as a discrete $\infty$-category. 
The procedure of taking a colimit gives a functor
	\[
	\mathsf{Fun}(\mathbb{Z}^{\mathrm{ds}}_{\ge 0}, \mathsf{Sp})
	\longrightarrow \mathsf{Sp}
	\]
given informally on objects by the formula $\{X_n\} \mapsto \bigoplus_{n\ge 0} X_n$.
We refer to this as the underlying spectrum of the graded spectrum $\{X_n\}$.

The homotopy theory
of graded spectra can be promoted to a symmetric monoidal $\infty$-category
via Day convolution, using addition to give $\mathbb{Z}_{\ge 0}^{\mathrm{ds}}$
a symmetric monoidal structure. With respect to this structure, we have
an equivalence
	\[
	\mathsf{Alg}_{\mathbb{E}_n}
	(\mathsf{Fun}(\mathbb{Z}^{\mathrm{ds}}_{\ge 0}, \mathsf{Sp}))
	\simeq
	\mathsf{Fun}^{\mathrm{lax-}\mathbb{E}_n}
	(\mathbb{Z}^{\mathrm{ds}}_{\ge 0}, \mathsf{Sp}),
	\]
where the latter $\infty$-category is the $\infty$-category of lax
$\mathbb{E}_n$-monoidal functors
$\mathbb{Z}^{\mathrm{ds}}_{\ge 0} \to \mathsf{Sp}$. We note that the
functor assigning to a graded spectrum its underlying spectrum can be
made symmetric
monoidal in an essentially unique way.

Now recall that the assignment $n \mapsto S^{2n}$ may be promoted
to an $\mathbb{E}_2$-monoidal functor
	\[
	\mathbb{Z}_{\ge 0}^{\mathrm{ds}} \longrightarrow
	\mathsf{Pic}(S^0).
	\]
(See, for example, \cite[Proposition 5.1.13]{rot}).
If $k\ge 0$ define $A_k \in \mathsf{Fun}(\mathbb{Z}_{\ge 0}^{\mathrm{ds}},
\mathsf{Sp})$ as the composite
	\[
	\xymatrix{
	\mathbb{Z}_{\ge 0}^{\mathrm{ds}} \ar[r]^{\cdot k} &
	\mathbb{Z}_{\ge 0} \ar[r] & \mathsf{Pic}(S^0)
	\ar[r] & \mathsf{Sp}
	}
	\]
If $k<0$, define $A_k$ as the composite
	\[
	\xymatrix{
	\mathbb{Z}_{\ge 0}^{\mathrm{ds}} \ar[r]^{\cdot (-k)} &
	\mathbb{Z}_{\ge 0} \ar[r] & \mathsf{Pic}(S^0)
	\ar[r]^{D}& \mathsf{Pic}(S^0)
	\ar[r] & \mathsf{Sp}
	}
	\]
where $D: \mathsf{Pic}(S^0) \to \mathsf{Pic}(S^0)$ denotes 
Spanier-Whitehead duality. In both cases, we see that $A_k$ is
$\mathbb{E}_2$-monoidal, being the composite of
$\mathbb{E}_2$-monoidal maps, and hence we may regard
$A_k$ as an $\mathbb{E}_2$-algebra object in
$\mathsf{Fun}(\mathbb{Z}^{\mathrm{ds}}_{\ge 0}, \mathsf{Sp})$. 
The underlying spectrum is readily checked to be 
$S^0[S^{2k}]$, so we have proven Lemma B.

We end by justifying the claim that nonnegatively graded $\mathbb{E}_n$-algebras
can be canonically augmented over their restriction to degree zero.
Denote by $\mathsf{Fun}(\mathbb{Z}^{\mathrm{ds}}_{\ge 0}, \mathsf{Sp})_{=0}$
the full subcategory of graded spectra which are concentrated in degree $0$. 
This is a localization of $\mathsf{Fun}(\mathbb{Z}^{\mathrm{ds}}_{\ge 0},
\mathsf{Sp})$ which is compatible with the symmetric monoidal structure
since, if $X \to Y$ is a map of nonnegatively graded spectra such that
$X_0 \to Y_0$ is an equivalence, and $Z$ is any other nonnegatively 
graded spectrum, then 
	\[
	(X \otimes Z)_0 \simeq X_0 \otimes Z_0 \to
	Y_0 \otimes Z_0 \simeq (Y \otimes Z)_0
	\]
is an equivalence. (Notice that this fails if we take $\mathbb{Z}$-graded
spectra, and, indeed, many $\mathbb{Z}$-graded algebra do not admit
augmentations). It now follows from \cite[2.2.1.9]{ha} that, if $A$ is
an $\mathbb{E}_n$-algebra in 
$\mathsf{Fun}(\mathbb{Z}^{\mathrm{ds}}_{\ge 0}, \mathsf{Sp})$, the counit
$A \to A_0$ is a map of graded $\mathbb{E}_n$-algebras. In particular, the
underlying $\mathbb{E}_n$-algebra of $A$ admits an $\mathbb{E}_n$-algebra
map to $A_0$.